\documentclass[12pt, reqno]{amsart}

\usepackage{graphicx}% Include figure files
\usepackage{amssymb}
\usepackage{amsmath}
\usepackage{bm}% bold math
\usepackage[left=2.5cm, right=2.5cm, top=3cm]{geometry}
\usepackage{hyperref}

\hypersetup{
    pdftitle={Boltzmann},
    pdfauthor={Macc},
    pdfmenubar=false,
    pdffitwindow=true,
    pdfstartview=FitH,
    colorlinks=true,
    linkcolor=blue,
    citecolor=green,
    urlcolor=cyan
}

\newtheorem{theorem}{Theorem}

\newtheorem{fact}{Fact}
\newtheorem{definition}{Definition}
\newtheorem{corollary}[theorem]{Corollary}
\newtheorem{proposition}[theorem]{Proposition}
\newtheorem{lemma}[theorem]{Lemma}

\newtheorem{axiom}{Axiom A.}

\newenvironment{boxer}
{\begin{center}
\begin{tabular}{|p{0.98\columnwidth}|}
\hline \\
}
{
\\ \\ \hline
\end{tabular}
\end{center}
}

\begin{document}

\title{Axiomatic Tests for the Boltzmann Distribution}% Force line breaks with \\
\thanks{
We wish to thank Riccardo Zecchina for helpful comments, Marco
Pirazzini and Giulio Principi for brilliant research assistance. Part of the material of this paper was first circulated in the IGIER working paper 593 of 2016.
}%

\author{Simone Cerreia-Vioglio}
\address{Universit\`a Bocconi}
\email{simone.cerreia@unibocconi.it}

\author{Fabio Maccheroni}
%\address{Department of Decision Sciences, Universit\`a Bocconi, via Roentgen 1, Milan 20136, Italy}
\address{Universit\`a Bocconi}
\email{fabio.maccheroni@unibocconi.it}

\author{Massimo Marinacci}
\address{Universit\`a Bocconi}
\email{massimo.marinacci@unibocconi.it}

\author{Aldo Rustichini}
\address{University of Minnesota}
\email{aldo.rustichini@gmail.com}

%\address{%
% Bocconi University, and University of Minnesota
%}%

\date{}

\begin{abstract}
The Boltzmann distribution describes a single parameter (temperature)
family of probability distributions over a state space; at any given temperature, the ratio of
probabilities of two states depends on their difference in energy. The same
family is known in other disciplines (economics, psychology, computer science)
with different names and interpretations.
Such widespread use in very diverse fields suggests a common
conceptual structure. We identify it on the basis of few natural axioms.
Checking whether observables satisfy these axioms is easy, so our characterization
provides a simple empirical test of the Boltzmannian modeling theories.
\end{abstract}

\maketitle

\section{Introduction}

According to the classic \emph{Boltzmann distribution} of statistical
mechanics, when the energy associated with some state $a$ of a system is
$E\left(  a\right)  $, then the frequency with which that state occurs in equilibrium is
proportional to%
\[
e^{-\frac{E\left(  a\right)  }{kt}}%
\]
where $t$ is the system absolute temperature and $k$ is the Boltzmann constant.

Under different interpretations and names (e.g., \emph{softmax} or
\emph{Multinomial Logit}), the Boltzmann distribution is widely used in many
fields of science, from physics to computer science, from economics to
psychology. For example, in economics the Multinomial Logit distribution is
the workhorse of discrete choice analysis. It gives the probability that an
agent with a utility function $V=-E$ selects an alternative $a$ when trying to
maximize $V$ but, say because of lack of information, makes mistakes in
evaluating the various alternatives. In this case the standard deviation of
mistakes is proportional to $t$.\footnote{See Train \cite{train} for a
textbook presentation. Later we will discuss another recent use of the
Multinomial Logit distribution in economics (Section \ref{sect:oia}).}
More recently, in econophysics, the Boltzmann distribution has been used to describe
market imperfections (with $E$ representing the bid-ask spreads of quotations)
and income distributions (with $E$ representing the amounts of money corresponding
to wealth levels).\footnote{See the letter of Kanazawa et al. \cite{KS},
 and  the colloquium of Yakovenko and Rosser \cite{YR}.} 

In this paper, we provide an axiomatic characterization of the Boltzmann
distribution based on observables. Specifically, we show that a family
$p=\left\{  p_{t}\right\}  $ of conditional distributions satisfies certain
properties if and only if there exists an energy function $E$ such that%
\[
p_{t}\left(  a\mid A\right)  =\dfrac{e^{-\frac{E\left(  a\right)  }{kt}}}%
{\sum_{b\in A}e^{-\frac{E\left(  b\right)  }{kt}}}%
\]
for all temperatures $t$ and all states $a$ in a collection $A$ of accessible
states. The function $E$ is unique up to an additive constant and can be
retrieved from data. Besides a common conceptual structure for this ubiquitous
distribution, our axiomatic analysis thus provides an empirical test for it.

\section{Notation}

We denote by $\mathcal{A}$ the collection of all finite subsets $A$ of a
universal system of states $X$, with $\left\vert X\right\vert \geq2$, and by
$p$ a \emph{random state function}%
\[%
\begin{array}
[c]{cccc}%
p: & \left(  0,\infty\right)  \times X\times\mathcal{A} & \rightarrow &
\mathbb{R}_{+}\\
& \left(  t,a,A\right)  & \mapsto & p_{t}\left(  a\mid A\right)
\end{array}
\]
that associates to a triplet $\left(  t,a,A\right)  $ the frequency
$p_{t}\left(  a\mid A\right)  $ of state $a\in X$, at\ temperature $t$ when
$A$ is the subsystem of accessible states.

Clearly,
\[
p_{t}\left(  B\mid A\right)  =\sum\nolimits_{b\in A}p_{t}\left(  b\mid
A\right)
\]
is the conditional frequency of some state in $B\subseteq A$. For a binary
subsystem, we just write $p_{t}\left(  a,b\right)  =p_{t}\left(  a\mid\left\{
a,b\right\}  \right)  $ for the frequency of a state $a$, with its odds
denoted by
\[
r_{t}\left(  a,b\right)  =\frac{p_{t}\left(  a,b\right)  }{p_{t}\left(
b,a\right)  }%
\]
%To ease notation, we indicate by $p_{t}\left(  B\right)  =p_{t}\left(  B\mid
%X\right)  $ the unconditional frequency of some state in $B\subseteq A$, and
%by just $p_{t}\left(  x\right)  $ the unconditional frequency of state $x\in
%X$.
Finally, $\delta_{a}$ is the point mass at $a\in X$, i.e., $\delta_{a}\left(
A\right)  =1$ if $a\in A$ and $\delta_{a}\left(  A\right)  =0$ otherwise.

\section{Axioms and results}

We consider the following axioms on a given random state function $p:\left(
0,\infty\right)  \times X\times\mathcal{A}\rightarrow\mathbb{R}_{+}$ that
describes the statistical behavior of the system.

We begin with positivity and conditioning axioms that require each section
$p_{t}$ of $p$ to be a conditional probability system (see Renyi \cite{Ren}
and Luce \cite{Luc}).

\begin{axiom}
[Positivity]\label{ax:pors}Given any $\left(  t,A\right)  \in\left(
0,\infty\right)  \times\mathcal{A}$,%
\[
\sum\nolimits_{a\in X}p_{t}\left(  a\mid A\right)  =1
\]
with $p_{t}\left(  a\mid A\right)  >0$ if and only if $a\in A$.
\end{axiom}

\begin{axiom}
[Conditioning]\label{ax:luce}Given any $\left(  t,A\right)  \in\left(
0,\infty\right)  \times\mathcal{A}$,%
\[
p_{t}\left(  b\mid A\right)  =p_{t}\left(  b\mid B\right)  p_{t}\left(  B\mid
A\right)
\]
\emph{\ }for all $B\subseteq A$ and all $b\in B$.
\end{axiom}

The next axiom requires the conditional probability systems $p_{t}$ to vary
continuously with temperature.

\begin{axiom}
[Continuity]\label{ax:cont}Given any $\left(  a,A\right)  \in X\times
\mathcal{A}$,%
\[
\lim_{t\rightarrow s}p_{t}\left(  a\mid A\right)
\]
exists for all $s\geq0$ and coincides with $p_{s}\left(  a\mid A\right)  $
when $s>0$.
\end{axiom}

Continuity guarantees, \emph{inter alia}, that as $t$ goes to $0$ a limit
probability $p_{0}\left(  a\mid A\right)  $ is defined for all $\left(
a,A\right)  \in X\times\mathcal{A}$. The following axiom requires the
consistency of freezing and positive temperature probabilities.

\begin{axiom}
[Consistency]\label{ax:cons}Given any $a,b\in X$,%
\[
p_{t}\left(  a,b\right)  >p_{t}\left(  b,a\right)  \implies p_{0}\left(
a,b\right)  >p_{0}\left(  b,a\right)
\]
for all $t>0$.
\end{axiom}

Next we postulate that, if at a zero temperature a binary subsystem is not
deterministically in either state, then both states are equally likely.

\begin{axiom}
[Zero Uniformity]\label{ax:unif}Given any $a,b\in X$,%
\[
p_{0}\left(  a,b\right)  \neq0,1\implies p_{0}\left(  a,b\right)  =1/2
\]

\end{axiom}

Our final axiom ties together the conditional distributions at different temperatures. It requires the dependence of odds from inverse temperature not
to be infinitely far from exponential. It is just a \textquotedblleft grain of
exponentiality\textquotedblright\ in the dependence of the system on time, that, as our next
theorem shows, develops into precisely an  exponential dependence of odds on inverse temperatures.

\begin{axiom}
[Boundedness]\label{ax:boun}Given any $a,b\in X$,%
\[
\sup_{t,s\in\left(  0,\infty\right)  }\left\vert r_{\frac{1}{t+s}}\left(
a,b\right)  -r_{\frac{1}{t}}\left(  a,b\right)  r_{\frac{1}{s}}\left(
a,b\right)  \right\vert <\infty
\]

\end{axiom}

We can now state our first result, in which we characterize the Boltzmann distribution.

\begin{theorem}
\label{thm:Boltzmann}A random state function $p:\left(  0,\infty
\right)  \times X\times\mathcal{A}\rightarrow\mathbb{R}_{+}$
satisfies A.\ref{ax:pors}--A.\ref{ax:boun} if and only if there exists a
function $E:X\rightarrow\mathbb{R}$ such that%
\begin{equation}
p_{t}\left(  a\mid A\right)  =\dfrac{e^{-\frac{E\left(  a\right)  }{kt}}}%
{\sum_{b\in A}e^{-\frac{E\left(  b\right)  }{kt}}}\delta_{a}\left(  A\right)
\label{eq:Boltzmann}%
\end{equation}
for all $\left(  t,a,A\right)  \in\left(  0,\infty\right)  \times
X\times\mathcal{A}$. Moreover, the function $E$ is unique up to an additive constant.
\end{theorem}

In view of this result, it is natural to say that a random state function $p$
is \emph{Boltzmannian} if it satisfies A.\ref{ax:pors}--A.\ref{ax:boun}. A
natural question is whether one can replace the thermal energy $kt$ with a
more general noise term $\kappa\left(  t\right)  $. To address this question,
we introduce a generic binary operation, concatenation, written $\oplus$, that has the usual sum $+$ as a
special case.

\begin{definition}
\label{def:conca}A \emph{concatenation }is a binary operation $\oplus$ on
$\mathbb{R}_{+}$\ which is associative, commutative, with identity element
$0$, and such that
\[
t>s\implies t\oplus v>s\oplus v\qquad\forall v\in\left(  0,\infty\right)
\]

\end{definition}

Besides the sum, other simple examples of concatenation are $t\oplus
s=t+s+\eta ts$ and $t\oplus s=\sqrt[\eta]{t^{\eta}+s^{\eta}}$ for some
$\eta\in\left(  0,\infty\right)  $.

The next axiom is the obvious extension of A.\ref{ax:boun} to a generic
concatenation. It continues to have a \textquotedblleft grain of
exponentiality\textquotedblright\ nature.

\begin{axiom}
[Weak Boundedness]\label{ax:webo}Given any $a,b\in X$,%
\[
\sup_{t,s\in\left(  0,\infty\right)  }\left\vert r_{\frac{1}{t\oplus s}%
}\left(  a,b\right)  -r_{\frac{1}{t}}\left(  a,b\right)  r_{\frac{1}{s}%
}\left(  a,b\right)  \right\vert <\infty
\]
for a continuous concatenation $\oplus$.
\end{axiom}

We can now generalize the earlier Boltzmannian result, which is the special
case of the theorem below when the concatenation $\oplus$ is the usual sum
$+$. A final notion is needed: $p$ is \emph{uniform} when $p_{t}\left(  a\mid
A\right)  =\delta_{a}\left(  A\right)  /\left\vert A\right\vert $ for all
$\left(  t,a,A\right)  \in\left(  0,\infty\right)  \times X\times\mathcal{A}$.

\begin{theorem}
\label{thm:Luce_on_steroids}A random state function $p:\left(
0,\infty\right)  \times X\times\mathcal{A}\rightarrow\mathbb{R}_{+}$
satisfies A.\ref{ax:pors}--A.\ref{ax:unif} and A.\ref{ax:webo} if and only if
there exist a function $E:X\rightarrow\mathbb{R}$ and an increasing bijection
$\kappa:\left(  0,\infty\right)  \rightarrow\left(  0,\infty\right)  $ such
that%
\begin{equation}
p_{t}\left(  a\mid A\right)  =\left\{
\begin{array}
[c]{lll}%
\dfrac{e^{-\frac{E\left(  a\right)  }{\kappa\left(  t\right)  }}}{\sum_{b\in
A}e^{-\frac{E\left(  b\right)  }{\kappa\left(  t\right)  }}}\bigskip &  & a\in
A\\
0 &  & a\notin A
\end{array}
\right.  \label{eq:softmax_txt}%
\end{equation}
for all $\left(  t,a,A\right)  \in\left(  0,\infty\right)  \times
X\times\mathcal{A}$.

In this case, $p$ is uniform if and only if $E$ is a constant function. When
$E$ is non-constant:

\begin{itemize}
\item[(i)] functions $\tilde{E}$ and $\tilde{\kappa}$ also represent $p$ as in
(\ref{eq:softmax_txt}) if and only if there exist $m>0$ and $q\in\mathbb{R}$
such that $\tilde{E}=mE+q$ and $\tilde{\kappa}=m\kappa$;
\end{itemize}

\begin{itemize}
\item[(ii)] the only concatenation $\oplus$ for which A.\ref{ax:webo} holds is%
\[
t\oplus s=\phi^{-1}\left[  \phi\left(  t\right)  +\phi\left(  s\right)
\right]  \qquad\forall t,s\in\left[  0,\infty\right)
\]
where $\phi:\left[  0,\infty\right)  \rightarrow\left[  0,\infty\right)  $ is
given by $\phi\left(  v\right)  =1/\kappa\left(  1/v\right)  $ for all $v>0$
and $\phi\left(  0\right)  =0$.
\end{itemize}
\end{theorem}

\section{Convex energy}

The physical question that the Boltzmann distribution addressed was: What is
the distribution of velocities in a gas at a certain temperature? The space of
states (velocities) is a convex set, and energy, which is proportional to square
speed, is a convex function. Analogously, in economics concave utility
functions play a fundamental role.

This motivates the next result that characterizes convex energy (so, concave utility).

\begin{proposition}
\label{prop:convexity}Let $X$ be a convex set and $p$ a Boltzmannian random
state function with energy $E$. The following conditions are equivalent:

\begin{enumerate}
\item[(i)] the function $E:X\rightarrow\mathbb{R}$ is convex;

\item[(ii)] there exists $t\in\left(  0,\infty\right)  $ such that%
\begin{equation}
p_{\alpha t}\left(  \alpha a+\left(  1-\alpha\right)  b,b\right)  \geq
p_{t}\left(  a,b\right)  \label{eq:baffo}%
\end{equation}
for all $a,b\in X$ and all $\alpha\in\left(  0,1\right)  $;

\item[(iii)] given any $s\in\left(  0,\infty\right)  $,%
\begin{equation}
p_{s}\left(  b\ \left\vert \ \frac{1}{\eta}A+\left(  1-\frac{1}{\eta}\right)
b\right.  \right)  \leq p_{\eta s}\left(  b\mid A\right)  \label{eq:baffo1}%
\end{equation}
for all $A\in\mathcal{A}$, all $b\in A$, and all $\eta>1$;

\item[(iv)] given any $s\in\left(  0,\infty\right)  $,%
\begin{equation}
p_{s}\left(  b\ \left\vert \ \frac{1}{\eta}A+\left(  1-\frac{1}{\eta}\right)
b\right.  \right)  \leq p_{\eta s}\left(  b\mid A\right)  \label{eq:barba}%
\end{equation}
for all $A\in\mathcal{A}$, all $b\in\arg\min_{a\in A}p_{\eta s}\left(  a\mid
A\right)  $, and all $\eta>1$.
\end{enumerate}
\end{proposition}

The stochastic choice interpretation of this result is based on the trade-off between noise (temperature) and states' distinguishability. By mixing states we make them closer, so less distinguishable and we augment the probability of making a mistake. To compensate such a mixing, according to inequality (\ref{eq:baffo}) it is more than sufficient to decrease noise proportionally. To illustrate, if $a$ is an optimal state with $\alpha=1/2$ the
inequality becomes%
\[
p_{\frac{t}{2}}\left(  \frac{1}{2}a+\frac{1}{2}b,b\right)  \geq p_{t}\left(
a,b\right)
\]
So, a even mixing is overcompensated by halving the noise. A similar
interpretation can be given to the other inequalities (\ref{eq:baffo1}) and
(\ref{eq:barba}).

\section{Additional remarks}

\subsection{Axioms' falsifiability}

As to the axioms' falsifiability, first observe that when $p$ is
\emph{uniform} the axioms hold, $E$ is constant, and $\kappa$ is undetermined.
The \emph{non-uniform} case depends on whether or not $p_{v}\left(
c,d\right)  =p_{v}\left(  d,c\right)  $ for all $v\in\left(  0,\infty\right)
$ and all $c\neq d$ in $X$. If this is the case, then Axiom A.\ref{ax:luce} is
violated.\footnote{If A.\ref{ax:luce} held, then, for all $A\in\mathcal{A}$
and all $c,d\in A$, it would follow%
\[
p_{v}\left(  c\mid A\right)  =p_{v}\left(  c,d\right)  p_{v}\left(  \left\{
c,d\right\}  \mid A\right)  =p_{v}\left(  d,c\right)  p_{v}\left(  \left\{
d,c\right\}  \mid A\right)  =p_{v}\left(  d\mid A\right)
\]
yielding uniformity of $p$.} Otherwise, we have the following result.

\begin{proposition}
\label{prop:bartsuka}Let $p$ be a random state function that satisfies
$p_{\bar{v}}\left(  \bar{c},\bar{d}\right)  >p_{\bar{v}}\left(  \bar{d}%
,\bar{c}\right)  $ for some $\bar{v}\in\left(  0,\infty\right)  $ and $\bar
{c},\bar{d}\in X$. If A.\ref{ax:pors}--A.\ref{ax:unif} are not violated, then
A.\ref{ax:webo} is satisfied if and only if representation
(\ref{eq:softmax_txt}) holds with
\[
\tilde{E}\left(  a\right)  =\ln r_{\bar{v}}\left(  \bar{c},a\right)
\qquad\text{and\qquad}\tilde{\kappa}\left(  t\right)  =\frac{\ln r_{\bar{v}%
}\left(  \bar{c},\bar{d}\right)  }{\ln r_{t}\left(  \bar{c},\bar{d}\right)  }%
\]
for all $\left(  t,a\right)  \in\left(  0,\infty\right)  \times X$.
\end{proposition}

Therefore, by Theorem \ref{thm:Luce_on_steroids}-(ii) the only concatenation
$\oplus$ for which A.\ref{ax:webo} holds corresponds to
\begin{equation}
\phi\left(  t\right)  =\frac{1}{\kappa\left(  1/t\right)  }=\frac{\ln
r_{1/t}\left(  \bar{c},\bar{d}\right)  }{\ln r_{\bar{v}}\left(  \bar{c}%
,\bar{d}\right)  }\qquad\forall t\in\left(  0,\infty\right)  \label{eq:conca}%
\end{equation}
In this way, observability of $\ln r_{1/t}\left(  \bar{c},\bar{d}\right)  $
qualifies the \textquotedblleft for some\textquotedblright\ clause of
A.\ref{ax:webo} and makes it falsifiable.

\subsection{Alternative axioms}

We can replace A.\ref{ax:cons} and A.\ref{ax:webo} in Theorem
\ref{thm:Luce_on_steroids} with the following two.

\begin{axiom}
[Monotonicity]\label{ax:mono}Given any $a,b\in X$, $\lim_{v\rightarrow\infty
}r_{v}\left(  a,b\right)  =1$; moreover,%
\[%
\begin{array}
[c]{l}%
r_{t}\left(  a,b\right)  >1\iff r_{s}\left(  a,b\right)  >r_{t}\left(
a,b\right)
\end{array}
\]
for all $s<t$\ in $\left(  0,\infty\right)  $.
\end{axiom}

\begin{axiom}
[Concatenation]\label{ax:conc}Given any $a,b,c,d\in X$,%
\[
r_{v}\left(  a,b\right)  >r_{t}\left(  a,b\right)  r_{s}\left(  a,b\right)
\implies r_{v}\left(  c,d\right)  >r_{t}\left(  c,d\right)  r_{s}\left(
c,d\right)
\]
\ for all $s,t,v\in\left(  0,\infty\right)  $\ such that $r_{v}\left(
a,b\right)  >1$\ and $r_{v}\left(  c,d\right)  >1.$
\end{axiom}

Next we establish the equivalence of these axioms with the earlier ones.

\begin{proposition}
\label{lem:operationsbis}Let $p:\left(  0,\infty\right)  \times X\times
\mathcal{A}\rightarrow\mathbb{R}_{+}$ be a random state function that
satisfies A.\ref{ax:pors}--A.\ref{ax:cont}, and A.\ref{ax:unif}. Then, $p$
satisfies A.\ref{ax:cons} and A.\ref{ax:webo} if and only if it satisfies
A.\ref{ax:mono} and A.\ref{ax:conc}.
\end{proposition}

Different sets of axioms for more general Multinomial Logit forms appear in the
subsequent papers of Saito \cite{Sa} and Cerreia-Vioglio et al. \cite{Ce}.

\subsection{Optimal information acquisition\label{sect:oia}}

In economics, the multinomial logit distribution has been used to formalize
versions of the discovered preference hypothesis, where the utility function
$V=-E$ is to be learned by an agent who confronts a cost $t$ of acquiring and
processing one unit of information. In particular, Matejka and McKay
\cite{matvejka2015rational} showed that the multinomial logit distribution
gives the optimal choice probability with which such an agent chooses an
alternative $a$ from a set $A$ of (a priori homogeneous) available
alternatives. Our axioms allow an analyst who controls $t$ to test
this theory.

In this economic setting, the concavity of the utility function is based on
the trade-off between decision time and alternatives' distinguishability. Now
inequality (\ref{eq:baffo}) says that, to compensate a mixing of alternatives
with a factor $\alpha$, which makes them less distinguishable, it is more than
sufficient to increase the decision time by a factor $1 / \alpha$.

\section{Proofs and related material}

A theorem of Aczel \cite{Acz} characterizes continuous concatenations.

\begin{theorem}
[Aczel]\label{thm:acz}A binary operation $\oplus$ on $\mathbb{R}_{+}$ is a
continuous concatenation if and only if there exists an increasing bijection
$f:\mathbb{R}_{+}\rightarrow\mathbb{R}_{+}$ such that
\[
t\oplus s=f^{-1}\left(  f\left(  t\right)  +f\left(  s\right)  \right)
\qquad\forall t,s\in\mathbb{R}_{+}%
\]
In this case, $f\left(  0\right)  =0$ and $f$ is strictly increasing and continuous.
\end{theorem}

The function $f$ is said be a generator for $\oplus$, which is then denoted by
$\oplus_{f}$.

\begin{lemma}
\label{lem:operations}If $p:\left(  0,\infty\right)  \times X\times
\mathcal{A}\rightarrow\mathbb{R}_{+}$ is a random state function that
satisfies A.\ref{ax:pors}, A.\ref{ax:cont}, A.\ref{ax:cons}, and
A.\ref{ax:unif}, then:

\begin{itemize}
\item[(i)] the relation defined on $X$ by $a\succsim b$ if and only if
$p_{0}\left(  a,b\right)  >0$ is such that%
\begin{align*}
\left.  a\succ b\right.   &  \iff p_{0}\left(  a,b\right)  >p_{0}\left(
b,a\right) \\
&  \iff p_{0}\left(  a,b\right)  =1\text{ and }p_{0}\left(  b,a\right)  =0\\
\left.  a\sim b\right.   &  \iff p_{0}\left(  a,b\right)  =p_{0}\left(
b,a\right) \\
&  \iff p_{0}\left(  a,b\right)  =p_{0}\left(  b,a\right)  \in\left\{
1,1/2\right\} \\
\left.  b\succ a\right.   &  \iff p_{0}\left(  a,b\right)  <p_{0}\left(
b,a\right) \\
&  \iff p_{0}\left(  a,b\right)  =0\text{ and }p_{0}\left(  b,a\right)  =1
\end{align*}

\item[(ii)] given any $a,b\in X$, the function $\varphi_{a,b}:\left(
0,\infty\right)  \rightarrow\left(  0,\infty\right)  $ defined by
\[
\varphi_{a,b}\left(  t\right)  =r_{1/t}\left(  a,b\right)  \qquad\forall
t\in\left(  0,\infty\right)
\]
is continuous and either diverges to $\infty$ as $t\rightarrow\infty$ (if
$a\succ b$) or is constantly equal to $1$ (if $a\sim b$) or vanishes as
$t\rightarrow\infty$ (if $b\succ a$).
\end{itemize}
\end{lemma}

\noindent\textbf{Proof} A.\ref{ax:pors} and A.\ref{ax:cont} imply that
$p_{0}\left(  \cdot\mid\left\{  a,b\right\}  \right)  $ is a probability
distribution (supported) on $\left\{  a,b\right\}  $, for all $a,b\in X$. The
proof is made pedantic by the fact that, if $a=b$, then $\left\{  a,b\right\}
=\left\{  a\right\}  =\left\{  b\right\}  $ and
\[
p_{0}\left(  a,b\right)  +p_{0}\left(  b,a\right)  =p_{0}\left(  a\mid\left\{
a\right\}  \right)  +p_{0}\left(  b\mid\left\{  b\right\}  \right)  =2
\]
else $a\neq b$ and
\[
p_{0}\left(  a,b\right)  +p_{0}\left(  b,a\right)  =1
\]

(i) By definition, $a\succ b$ iff $a\succsim b$ and not $b\succsim a$, that
is, $p_{0}\left(  a,b\right)  >0$ and $p_{0}\left(  b,a\right)  \leq0$.

\begin{itemize}
\item Assume $a\succ b$, then $p_{0}\left(  a,b\right)  >0$ and $p_{0}\left(
b,a\right)  \leq0$ imply $p_{0}\left(  a,b\right)  >p_{0}\left(  b,a\right)  $.
\end{itemize}

\begin{itemize}
\item Assume $p_{0}\left(  a,b\right)  >p_{0}\left(  b,a\right)  $. This is
impossible if $a=b$, therefore $a\neq b$ and $p_{0}\left(  a,b\right)  >0$. If
it held $p_{0}\left(  b,a\right)  >0$, then $p_{0}\left(  a,b\right)
+p_{0}\left(  b,a\right)  =1$ would imply $p_{0}\left(  a,b\right)
,p_{0}\left(  b,a\right)  \in\left(  0,1\right)  $, and A.\ref{ax:unif} would
yield $p_{0}\left(  a,b\right)  =1/2=p_{0}\left(  b,a\right)  $, a
contradiction. Then it must be $p_{0}\left(  b,a\right)  =0$ and $p_{0}\left(
a,b\right)  =p_{0}\left(  a,b\right)  +p_{0}\left(  b,a\right)  =1$.
\end{itemize}

\begin{itemize}
\item Assume $p_{0}\left(  a,b\right)  =1$ and $p_{0}\left(  b,a\right)  =0$,
then $p_{0}\left(  a,b\right)  =1$ and $p_{0}\left(  b,a\right)  \leq0$, and
$a\succ b$.
\end{itemize}

By definition, $a\sim b$ iff $a\succsim b$ and also $b\succsim a$, that is,
$p_{0}\left(  a,b\right)  >0$ and $p_{0}\left(  b,a\right)  >0$.

\begin{itemize}
\item Assume $a\sim b$. If $a=b$, then $p_{0}\left(  a,b\right)
=1=p_{0}\left(  b,a\right)  $. Else $a\neq b$, $p_{0}\left(  a,b\right)
,p_{0}\left(  b,a\right)  >0$, and $p_{0}\left(  a,b\right)  +p_{0}\left(
b,a\right)  =1$, then $p_{0}\left(  a,b\right)  ,p_{0}\left(  b,a\right)
\in\left(  0,1\right)  $, and A.\ref{ax:unif} yields $p_{0}\left(  a,b\right)
=1/2=p_{0}\left(  b,a\right)  $.
\end{itemize}

\begin{itemize}
\item Assume $p_{0}\left(  a,b\right)  =p_{0}\left(  b,a\right)  $. If $a=b$,
then $p_{0}\left(  a,b\right)  =p_{0}\left(  b,a\right)  =1$. Else $a\neq b$,
and $p_{0}\left(  a,b\right)  +p_{0}\left(  b,a\right)  =1$, then
$2p_{0}\left(  a,b\right)  =1$ and $2p_{0}\left(  b,a\right)  =1$, that is,
$p_{0}\left(  a,b\right)  =p_{0}\left(  b,a\right)  =1/2$.
\end{itemize}

\begin{itemize}
\item Assume $p_{0}\left(  a,b\right)  =p_{0}\left(  b,a\right)  \in\left\{
1,1/2\right\}  $, then $p_{0}\left(  a,b\right)  ,p_{0}\left(  b,a\right)
>0$, and $a\sim b$.
\end{itemize}

The case $b\succ a$ follows from the case $a\succ b$ exchanging the roles of
the states.

(ii) Given any $t\in\left(  0,\infty\right)  $, $\varphi_{a,b}\left(
t\right)  =r_{1/t}\left(  a,b\right)  =p_{1/t}\left(  a,b\right)
/p_{1/t}\left(  b,a\right)  \in\left(  0,\infty\right)  $ for all $a,b\in X$
because $p_{1/t}\left(  \cdot\mid\left\{  a,b\right\}  \right)  $ is a
positive probability distribution on $\left\{  a,b\right\}  $, thus
$\varphi_{a,b}:\left(  0,\infty\right)  \rightarrow\left(  0,\infty\right)  $
is well defined. Moreover, by A.\ref{ax:cont}, $\varphi_{a,b}$ is also
continuous on $\left(  0,\infty\right)  $.

\begin{itemize}
\item If $a\succ b$, then $p_{0}\left(  a,b\right)  =1$ and $p_{0}\left(
b,a\right)  =0$, so $a\neq b$ and%
\[
\lim_{t\rightarrow\infty}\varphi_{a,b}\left(  t\right)  =\lim_{t\rightarrow
\infty}\frac{p_{1/t}\left(  a,b\right)  }{p_{1/t}\left(  b,a\right)  }%
=\lim_{t\rightarrow\infty}\frac{1-p_{1/t}\left(  b,a\right)  }{p_{1/t}\left(
b,a\right)  }=\infty
\]
hence $\varphi_{a,b}$ diverges at $\infty$ as $t\rightarrow\infty$.
\end{itemize}

For later reference, note that so far A.\ref{ax:cons} has not been used.

\begin{itemize}
\item If $a\sim b$, and \emph{per contra} $\varphi_{a,b}\left(  t\right)
\neq1$ for some $t\in\left(  0,\infty\right)  $, then

\begin{itemize}
\item[$\circ$] either $\varphi_{a,b}\left(  t\right)  >1$, thus $p_{1/t}%
\left(  a,b\right)  >p_{1/t}\left(  b,a\right)  $ and, by A.\ref{ax:cons},
$p_{0}\left(  a,b\right)  >p_{0}\left(  b,a\right)  $, contradicting $a\sim b$,

\item[$\circ$] or $\varphi_{a,b}\left(  t\right)  <1$, thus $p_{1/t}\left(
a,b\right)  <p_{1/t}\left(  b,a\right)  $ and, by A.\ref{ax:cons},
$p_{0}\left(  a,b\right)  <p_{0}\left(  b,a\right)  $, contradicting $a\sim b$,
\end{itemize}

in conclusion, $\varphi_{a,b}\left(  t\right)  =1$ for all $t\in\left(
0,\infty\right)  $.

\item If $b\succ a$, the thesis follows because $\varphi_{a,b}=1/\varphi
_{b,a}$.\hfill$\blacksquare\bigskip$
\end{itemize}

\noindent\textbf{Proof of Theorem \ref{thm:Luce_on_steroids}} Let $p$ be
a random state function that satisfies A.\ref{ax:pors}--A.\ref{ax:unif} and
A.\ref{ax:webo}. As in Lemma \ref{lem:operations}, define, for all $a,b\in X$,%
\[
\varphi_{a,b}\left(  t\right)  =r_{1/t}\left(  a,b\right)  \qquad\forall
t\in\left(  0,\infty\right)
\]
Also let $f:\left[  0,\infty\right)  \rightarrow\left[  0,\infty\right)  $ be
a generator of a continuous concatenation $\oplus=\oplus_{f}$ for which
A.\ref{ax:webo} holds. Set $g=f^{-1}$. By Theorem \ref{thm:acz}, $g:\left[
0,\infty\right)  \rightarrow\left[  0,\infty\right)  $ is a continuous and
strictly increasing bijection such that $g\left(  0\right)  =0$.

Next we show that, given any $a,b\in X$,%
\begin{equation}
\varphi_{a,b}\left(  g\left(  t+s\right)  \right)  =\varphi_{a,b}\left(
g\left(  t\right)  \right)  \varphi_{a,b}\left(  g\left(  s\right)  \right)
\qquad\forall t,s\in\left(  0,\infty\right)  \label{eq:phi-mult}%
\end{equation}
Three cases have to be considered, depending on whether $a\succ b$, $a\sim b$,
or $b\succ a$ according to the relation $\succsim$ defined in Lemma
\ref{lem:operations}.

\begin{itemize}
\item If $a\succ b$, then $\varphi_{a,b}$ is unbounded above and so is
$\varphi_{a,b}\circ g:\left(  0,\infty\right)  \rightarrow\left(
0,\infty\right)  $. Moreover, by A.\ref{ax:webo}, there exists $M>0$ such
that, for all $t,s\in\left(  0,\infty\right)  $%
\begin{align*}
\left\vert r_{\frac{1}{t\oplus s}}\left(  a,b\right)  -r_{\frac{1}{t}}\left(
a,b\right)  r_{\frac{1}{s}}\left(  a,b\right)  \right\vert  &  <M\\
\left\vert r_{\frac{1}{g\left(  g^{-1}\left(  t\right)  +g^{-1}\left(
s\right)  \right)  }}\left(  a,b\right)  -r_{\frac{1}{t}}\left(  a,b\right)
r_{\frac{1}{s}}\left(  a,b\right)  \right\vert  &  <M
\end{align*}
hence, for all $t^{\prime},s^{\prime}\in\left(  0,\infty\right)  $, choosing
$t=g\left(  t^{\prime}\right)  $ and $s=g\left(  s^{\prime}\right)  $, we
have
\begin{align*}
\left\vert r_{\frac{1}{g\left(  g^{-1}\left(  g\left(  t^{\prime}\right)
\right)  +g^{-1}\left(  g\left(  s^{\prime}\right)  \right)  \right)  }%
}\left(  a,b\right)  -r_{\frac{1}{g\left(  t^{\prime}\right)  }}\left(
a,b\right)  r_{\frac{1}{g\left(  s^{\prime}\right)  }}\left(  a,b\right)
\right\vert <M  & \\
\left\vert r_{\frac{1}{g\left(  t^{\prime}+s^{\prime}\right)  }}\left(
a,b\right)  -r_{\frac{1}{g\left(  t^{\prime}\right)  }}\left(  a,b\right)
r_{\frac{1}{g\left(  s^{\prime}\right)  }}\left(  a,b\right)  \right\vert <M
& \\
\left\vert \varphi_{a,b}\left(  g\left(  t^{\prime}+s^{\prime}\right)
\right)  -\varphi_{a,b}\left(  g\left(  t^{\prime}\right)  \right)
\varphi_{a,b}\left(  g\left(  s^{\prime}\right)  \right)  \right\vert <M  &
\end{align*}
But $\left(  0,\infty\right)  $ is a semigroup with respect to usual addition
and $\varphi_{a,b}\circ g$ is unbounded above. Therefore, Theorem 1 of Baker
\cite{Bak} implies that (\ref{eq:phi-mult}) holds.

\item If $a\sim b$, then $\varphi_{a,b}\left(  t\right)  =1$ for all
$t\in\left(  0,\infty\right)  $ and (\ref{eq:phi-mult}) holds.

\item Else, $b\succ a$ and, as the first point shows,%
\[
\varphi_{b,a}\left(  g\left(  t+s\right)  \right)  =\varphi_{b,a}\left(
g\left(  t\right)  \right)  \varphi_{b,a}\left(  g\left(  s\right)  \right)
\]
for all $t,s\in\left(  0,\infty\right)  $, but then
\begin{align*}
\varphi_{a,b}\left(  g\left(  t+s\right)  \right)   &  =\frac{1}{\varphi
_{b,a}\left(  g\left(  t+s\right)  \right)  }\\
&  =\frac{1}{\varphi_{b,a}\left(  g\left(  t\right)  \right)  \varphi
_{b,a}\left(  g\left(  s\right)  \right)  }\\
&  =\varphi_{a,b}\left(  g\left(  t\right)  \right)  \varphi_{a,b}\left(
g\left(  s\right)  \right)
\end{align*}
and (\ref{eq:phi-mult}) holds again.
\end{itemize}

Summing up, the functional equation (\ref{eq:phi-mult}) holds for all $a,b\in
X$. Continuity of $\varphi_{a,b}\circ g$, its strict positivity, and
(\ref{eq:phi-mult}), imply that
\[
\varphi_{a,b}\left(  g\left(  t\right)  \right)  =e^{v\left(  a,b\right)
t}\qquad\forall t\in\left(  0,\infty\right)
\]
for a unique $v\left(  a,b\right)  \in\mathbb{R}$ (see, e.g., Theorem 2.1.2.1
of Aczel \cite{Ac2}). It follows that $\varphi_{a,b}\left(  s\right)
=\varphi_{a,b}\left(  g\left(  f\left(  s\right)  \right)  \right)
=e^{v\left(  a,b\right)  f\left(  s\right)  }$ for all $s\in\left(
0,\infty\right)  $.

Now fix some $a^{\ast}\in X$ and define $E:X\rightarrow\mathbb{R}$ by
$E\left(  x\right)  =-v\left(  x,a^{\ast}\right)  $ for all $x\in X$. Given
any $t\in\left(  0,\infty\right)  $ and any $x,y\in X$, by A.\ref{ax:pors},
A.\ref{ax:luce}, and Theorem 2 of Luce \cite{Luc}, it follows that%
\begin{align*}
\varphi_{x,y}\left(  t\right)   &  =r_{1/t}\left(  x,y\right)  =r_{1/t}\left(
x,a^{\ast}\right)  r_{1/t}\left(  a^{\ast},y\right)  =\frac{r_{1/t}\left(
x,a^{\ast}\right)  }{r_{1/t}\left(  y,a^{\ast}\right)  }\\
&  =\frac{\varphi_{x,a^{\ast}}\left(  t\right)  }{\varphi_{y,a^{\ast}}\left(
t\right)  }=\frac{e^{v\left(  x,a^{\ast}\right)  f\left(  t\right)  }%
}{e^{v\left(  y,a^{\ast}\right)  f\left(  t\right)  }}=\frac{e^{-E\left(
x\right)  f\left(  t\right)  }}{e^{-E\left(  y\right)  f\left(  t\right)  }}%
\end{align*}
By Theorem 3 of Luce \cite{Luc}, for every $t\in\left(  0,\infty\right)  $,
$A\in\mathcal{A}$, and $a\in A$, arbitrarily choosing $c^{\ast}\in A$,%
\begin{align*}
p_{t}\left(  a\mid A\right)   &  =\frac{r_{t}\left(  a,c^{\ast}\right)  }%
{\sum_{b\in A}r_{t}\left(  b,c^{\ast}\right)  }=\frac{\varphi_{a,c^{\ast}%
}\left(  1/t\right)  }{\sum_{b\in A}\varphi_{b,c^{\ast}}\left(  1/t\right)
}\\
&  =\frac{\frac{e^{-E\left(  a\right)  f\left(  1/t\right)  }}{e^{-E\left(
c^{\ast}\right)  f\left(  1/t\right)  }}}{\sum_{b\in A}\frac{e^{-E\left(
b\right)  f\left(  1/t\right)  }}{e^{-E\left(  c^{\ast}\right)  f\left(
1/t\right)  }}}=\dfrac{e^{-f\left(  \frac{1}{t}\right)  E\left(  a\right)  }%
}{\sum_{b\in A}e^{-f\left(  \frac{1}{t}\right)  E\left(  b\right)  }}%
\end{align*}
and (\ref{eq:softmax_txt}) holds for $\kappa\left(  t\right)  =1/f\left(
1/t\right)  $ (because $p_{t}\left(  a\mid A\right)  =0$ for $a\notin A$ by
A.\ref{ax:pors}).

\begin{boxer}
\textbf{NB 1} So far we have shown that: If $p$ is random state function that
satisfies A.\ref{ax:pors}--A.\ref{ax:unif} and A.\ref{ax:webo} (with respect
to $\oplus_{f}$); then, setting $\kappa\left(  t\right)  =1/f\left(
1/t\right)  $ for all $t\in\left(  0,\infty\right)  $, there exists
$E:X\rightarrow\mathbb{R}$ such that%
\[
p_{t}\left(  a\mid A\right)  =\dfrac{e^{-\frac{E\left(  a\right)  }%
{\kappa\left(  t\right)  }}}{\sum_{b\in A}e^{-\frac{E\left(  b\right)
}{\kappa\left(  t\right)  }}}\delta_{a}\left(  A\right)
\]
for all $\left(  t,a,A\right)  \in\left(  0,\infty\right)  \times
X\times\mathcal{A}$.

Moreover, since $f_{|\left(  0,\infty\right)  }$ is a continuous and strictly
increasing bijection from $\left(  0,\infty\right)  $ to $\left(
0,\infty\right)  $, and $s\mapsto1/s$ is a continuous and strictly decreasing
bijection from $\left(  0,\infty\right)  $ to $\left(  0,\infty\right)  $,
then $\kappa:t\mapsto1/f\left(  1/t\right)  $ a continuous and strictly
increasing bijection from $\left(  0,\infty\right)  $ to $\left(
0,\infty\right)  $.
\end{boxer}

This proves the \textquotedblleft only if\textquotedblright\ part of the
statement.\bigskip

As to the \textquotedblleft if\textquotedblright\ part, assume that
(\ref{eq:softmax_txt}) holds. It is routine to check that $p$ satisfies
A.\ref{ax:pors}--A.\ref{ax:unif}. To prove that also A.\ref{ax:webo} holds,
define $\phi:\left[  0,\infty\right)  \rightarrow\left[  0,\infty\right)  $ by
setting $\phi\left(  v\right)  =1/\kappa\left(  1/v\right)  $ for all $v>0$,
and $\phi\left(  0\right)  =0$. Since $k$ is an increasing bijection from
$\left(  0,\infty\right)  $ to $\left(  0,\infty\right)  $, so is
$\phi_{|\left(  0,\infty\right)  }$. But then $\phi:\left[  0,\infty\right)
\rightarrow\left[  0,\infty\right)  $ is an increasing bijection too. Then
\begin{equation}
t\oplus s=\phi^{-1}\left[  \phi\left(  t\right)  +\phi\left(  s\right)
\right]  \qquad\forall t,s\in\left[  0,\infty\right)  \label{eq:oplus}%
\end{equation}
is a (well defined) binary operation on $\mathbb{R}_{+}$. Theorem
\ref{thm:acz} guarantees that $\oplus=\oplus_{\phi}$ is indeed a continuous
concatenation. With this, given $a,b\in X$, for all $t,s\in\left(
0,\infty\right)  $%
\[
r_{\frac{1}{t\oplus s}}\left(  a,b\right)  =e^{-\frac{1}{\kappa\left(
\frac{1}{t\oplus s}\right)  }\left[  E\left(  a\right)  -E\left(  b\right)
\right]  }%
\]
but, by (\ref{eq:oplus}), $t\oplus s>0$, hence, by definition of $\phi$,
\begin{align*}
r_{\frac{1}{t\oplus s}}\left(  a,b\right)   &  =e^{-\phi\left(  t\oplus
s\right)  \left[  E\left(  a\right)  -E\left(  b\right)  \right]
}=e^{-\left(  \phi\left(  t\right)  +\phi\left(  s\right)  \right)  \left[
E\left(  a\right)  -E\left(  b\right)  \right]  }\\
&  =e^{-\phi\left(  t\right)  \left[  E\left(  a\right)  -E\left(  b\right)
\right]  }e^{-\phi\left(  s\right)  \left[  E\left(  a\right)  -E\left(
b\right)  \right]  }\\
&  =e^{-\frac{1}{\kappa\left(  1/t\right)  }\left[  E\left(  a\right)
-E\left(  b\right)  \right]  }e^{-\frac{1}{\kappa\left(  1/s\right)  }\left[
E\left(  a\right)  -E\left(  b\right)  \right]  }\\
& =r_{1/t}\left(  a,b\right)
r_{1/s}\left(  a,b\right)
\end{align*}
A fortiori, A.\ref{ax:webo} holds, with respect to $\oplus_{\phi}$, where
$\phi\left(  v\right)  =1/\kappa\left(  1/v\right)  $ for all $v>0$, and
$\phi\left(  0\right)  =0$. Actually, we proved a stronger fact:

\begin{boxer}
\textbf{NB 2} Given a function $E:X\rightarrow\mathbb{R}$ and an increasing
bijection $\kappa:\left(  0,\infty\right)  \rightarrow\left(  0,\infty\right)
$, the function defined by
\[
p_{t}\left(  a\mid A\right)  =\dfrac{e^{-\frac{E\left(  a\right)  }%
{\kappa\left(  t\right)  }}}{\sum_{b\in A}e^{-\frac{E\left(  b\right)
}{\kappa\left(  t\right)  }}}\delta_{a}\left(  A\right)
\]
for all $\left(  t,a,A\right)  \in\left(  0,\infty\right)  \times
X\times\mathcal{A}$ is a random state function that satisfies A.\ref{ax:pors}%
--A.\ref{ax:unif} and A.\ref{ax:webo} (with respect to $\oplus_{\phi}$, where
$\phi\left(  v\right)  =1/\kappa\left(  1/v\right)  $ for all $v>0$, and
$\phi\left(  0\right)  =0$).
\end{boxer}

This concludes the proof of the first part of the statement.\bigskip

Now assume that (\ref{eq:softmax_txt}) holds for a function $E:X\rightarrow
\mathbb{R}$ and an increasing bijection $\kappa:\left(  0,\infty\right)
\rightarrow\left(  0,\infty\right)  $. Note that
\begin{align*}
r_{t}\left(  a,b\right)   &  =\exp\left(  -\frac{1}{\kappa\left(  t\right)
}\left[  E\left(  a\right)  -E\left(  b\right)  \right]  \right) \\
\frac{1}{\kappa\left(  t\right)  }\left[  E\left(  a\right)  -E\left(
b\right)  \right]   &  =-\ln r_{t}\left(  a,b\right)  =\ln r_{t}\left(
b,a\right)
\end{align*}
for all $\left(  t,a,b\right)  \in\left(  0,\infty\right)  \times X^{2}$.

If $p$ is uniform, then $\ln r_{t}\left(  b,a\right)  =0$ for all $\left(
t,a,b\right)  \in\left(  0,\infty\right)  \times X^{2}$, and strict positivity
of $\kappa$ implies $E$ is constant. The converse follows immediately from
(\ref{eq:softmax_txt}).

Else, $E$ is not constant. Let $\tilde{E}:X\rightarrow\mathbb{R}$ and
$\tilde{\kappa}:\left(  0,\infty\right)  \rightarrow\left(  0,\infty\right)  $
also represent $p$ as in (\ref{eq:softmax_txt}), then
\[
\frac{1}{\kappa\left(  t\right)  }\left[  E\left(  a\right)  -E\left(
b\right)  \right]  =\ln r_{t}\left(  b,a\right)  =\frac{1}{\tilde{\kappa
}\left(  t\right)  }\left[  \tilde{E}\left(  a\right)  -\tilde{E}\left(
b\right)  \right]
\]
for all $\left(  t,a,b\right)  \in\left(  0,\infty\right)  \times X^{2}$.
Arbitrarily choose $\left(  t^{\ast},a^{\ast},b^{\ast}\right)  \in\left(
0,\infty\right)  \times X^{2}$ such that $E\left(  a^{\ast}\right)  >E\left(
b^{\ast}\right)  $. Then:

(i) For all $a\in A$,%
\[
\frac{1}{\kappa\left(  t^{\ast}\right)  }\left[  E\left(  a\right)  -E\left(
b^{\ast}\right)  \right]  =\ln r_{t^{\ast}}\left(  b^{\ast},a\right)
=\frac{1}{\tilde{\kappa}\left(  t^{\ast}\right)  }\left[  \tilde{E}\left(
a\right)  -\tilde{E}\left(  b^{\ast}\right)  \right]
\]
hence
\[
\tilde{E}\left(  a\right)  =\underset{m^{\ast}}{\underbrace{\frac
{\tilde{\kappa}\left(  t^{\ast}\right)  }{\kappa\left(  t^{\ast}\right)  }}%
}E\left(  a\right)  +~\underset{q^{\ast}}{\underbrace{\tilde{E}\left(
b^{\ast}\right)  -\frac{\tilde{\kappa}\left(  t^{\ast}\right)  }{\kappa\left(
t^{\ast}\right)  }E\left(  b^{\ast}\right)  }}%
\]
and, for all $t\in\left(  0,\infty\right)  $,%
\begin{align*}
\frac{1}{\kappa\left(  t\right)  }\left[  E\left(  a^{\ast}\right)  -E\left(
b^{\ast}\right)  \right]   &  =\ln r_{t}\left(  b^{\ast},a^{\ast}\right)
=\frac{1}{\tilde{\kappa}\left(  t\right)  }\left[  \tilde{E}\left(  a^{\ast
}\right)  -\tilde{E}\left(  b^{\ast}\right)  \right]  \\
\frac{1}{\kappa\left(  t\right)  }\left[  E\left(  a^{\ast}\right)  -E\left(
b^{\ast}\right)  \right]   &  =\frac{1}{\tilde{\kappa}\left(  t\right)
}\left[  m^{\ast}E\left(  a^{\ast}\right)  -m^{\ast}E\left(  b^{\ast}\right)
\right]  \\
\tilde{\kappa}\left(  t\right)   &  =m^{\ast}\kappa\left(  t\right)
\end{align*}
thus there exist $m>0$ and $q\in\mathbb{R}$ such that $\tilde{E}=mE+q$ and
$\tilde{\kappa}=m\kappa$. This proves the \textquotedblleft only
if\textquotedblright\ part of point (i). The \textquotedblleft if\textquotedblright\  part is trivial.

(ii) By NB 2, under (\ref{eq:softmax_txt}), the binary operation defined by%
\[
t\oplus_{\phi}s=\phi^{-1}\left[  \phi\left(  t\right)  +\phi\left(  s\right)
\right]  \qquad\forall t,s\in\left[  0,\infty\right)
\]
where $\phi:\left[  0,\infty\right)  \rightarrow\left[  0,\infty\right)  $ is
given by $\phi\left(  v\right)  =1/\kappa\left(  1/v\right)  $ for all $v>0$
and $\phi\left(  0\right)  =0$, is a concatenation for which A.\ref{ax:webo}
holds. By NB 1, if $\oplus=\oplus_{f}$ is a concatenation for which
A.\ref{ax:webo} holds, then setting $\tilde{\kappa}\left(  t\right)
=1/f\left(  1/t\right)  $ for all $t\in\left(  0,\infty\right)  $, there
exists $\tilde{E}:X\rightarrow\mathbb{R}$ such that%
\[
p_{t}\left(  a\mid A\right)  =\dfrac{e^{-\frac{\tilde{E}\left(  a\right)
}{\tilde{\kappa}\left(  t\right)  }}}{\sum_{b\in A}e^{-\frac{\tilde{E}\left(
b\right)  }{\tilde{\kappa}\left(  t\right)  }}}\delta_{a}\left(  A\right)
\]
for all $\left(  t,a,A\right)  \in\left(  0,\infty\right)  \times
X\times\mathcal{A}$. By point (i), there exist $m>0$ and $q\in\mathbb{R}$ such
that $\tilde{E}=mE+q$ and $\tilde{\kappa}=m\kappa$; therefore, for all
$t\in\left(  0,\infty\right)  $,%
\[
\frac{1}{f\left(  1/t\right)  }=\tilde{\kappa}\left(  t\right)  =m\kappa
\left(  t\right)  =\frac{m}{\phi\left(  1/t\right)  }%
\]
hence $f=\phi/m$ on $\left(  0,\infty\right)  $, and $f\left(  0\right)
=0=\phi\left(  0\right)  /m$ by Theorem \ref{thm:acz}. Finally, $f=\phi/m$
implies $\oplus_{f}=\oplus_{\phi}$, concluding the proof of (ii).\hfill
$\blacksquare\bigskip$

\noindent\textbf{Proof of Theorem \ref{thm:Boltzmann} }If $p$ is a random state
function that satisfies A.\ref{ax:pors}--A.\ref{ax:unif} and A.\ref{ax:boun},
then it also satisfies A.\ref{ax:webo} with respect to $\oplus_{f}$\ where
$f\left(  t\right)  =t/k$ and $k$ is the Boltzmann constant. By NB 1 of the
previous proof, setting $\kappa\left(  t\right)  =1/f\left(  1/t\right)  $ for
all $t\in\left(  0,\infty\right)  $, it follows $\kappa\left(  t\right)  =kt$
and there exists $E:X\rightarrow\mathbb{R}$ such that%
\[
p_{t}\left(  a\mid A\right)  =\dfrac{e^{-\frac{E\left(  a\right)  }{kt}}}%
{\sum_{b\in A}e^{-\frac{E\left(  b\right)  }{kt}}}\delta_{a}\left(  A\right)
\]
for all $\left(  t,a,A\right)  \in\left(  0,\infty\right)  \times
X\times\mathcal{A}$. The converse is routine.

As to uniqueness of the representation, by point (i) of Theorem
\ref{thm:Luce_on_steroids}, if $\tilde{E}:X\rightarrow\mathbb{R}$, and
$\tilde{\kappa}\left(  t\right)  =kt$, also represent $p$ as in
(\ref{eq:Boltzmann}), then there exist $m>0$ and $q\in\mathbb{R}$ such that
$\tilde{E}=mE+q$ and $\tilde{\kappa}=m\kappa$, but this means $kt=mkt$ for all
$t>0$, that is $m=1$. Again, the converse is routine.\hfill$\blacksquare
\bigskip$

\noindent\textbf{Proof of Proposition \ref{prop:convexity}} \emph{(ii) is
equivalent to (i).} There exists $t\in\left(  0,\infty\right)  $ such that
(\ref{eq:baffo}) holds if and only if%
\begin{align*}
&  \left.  \exists t\in\left(  0,\infty\right)  :p_{\alpha t}\left(  \alpha
a+\left(  1-\alpha\right)  b,b\right)  \geq p_{t}\left(  a,b\right)  \right.
\\
&  \iff\left.  \exists t\in\left(  0,\infty\right)  :r_{\alpha t}\left(
\alpha a+\left(  1-\alpha\right)  b,b\right)  \geq r_{t}\left(  a,b\right)
\right. \\
&  \iff\left.  \exists t\in\left(  0,\infty\right)  :r_{\alpha t}\left(
b,\alpha a+\left(  1-\alpha\right)  b\right)  \leq r_{t}\left(  b,a\right)
\right. \\
&  \iff\left.  \exists t\in\left(  0,\infty\right)  :e^{\frac{1}{k\alpha
t}\left[  E\left(  \alpha a+\left(  1-\alpha\right)  b\right)  -E\left(
b\right)  \right]  }\leq e^{\frac{1}{kt}\left[  E\left(  a\right)  -E\left(
b\right)  \right]  }\right. \\
&  \iff E\left(  \alpha a+\left(  1-\alpha\right)  b\right)  \leq\alpha
E\left(  a\right)  +\left(  1-\alpha\right)  E\left(  b\right)
\end{align*}
for all $\left(  a,b,\alpha\right)  \in X\times X\times\left(  0,1\right)  $.

\emph{(i) implies (iii).} Given any $s\in\left(  0,\infty\right)  $,
$A\in\mathcal{A}$, $b\in A$,
and $\eta>1$,
\[
p_{s}\left(  b\ \left\vert \ \frac{1}{\eta}A+\left(  1-\frac{1}{\eta}\right)
b\right.  \right)  =\frac{1}{\sum_{a\in A}e^{-\frac{1}{ks}\left[  E\left(
\frac{1}{\eta}a+\left(  1-\frac{1}{\eta}\right)  b\right)  -E\left(  b\right)
\right]  }}%
\]
but convexity of $E$ implies $E\left(  \left(  1/\eta\right)  a+\left(
1-\left(  1/\eta\right)  \right)  b\right)  -E\left(  b\right)  \leq\left(
1/\eta\right)  \left(  E\left(  a\right)  -E\left(  b\right)  \right)  $
hence
\[
-\frac{1}{ks}\left[  E\left(  \left(  1/\eta\right)  a+\left(  1-\left(
1/\eta\right)  \right)  b\right)  -E\left(  b\right)  \right]  \geq-\frac
{1}{k\eta s}\left[  E\left(  a\right)  -E\left(  b\right)  \right]
\]
for all $a\in A$,
%	\begin{align*}
%	e^{-\frac{1}{ks}\left[  E\left(  \frac{1}{\eta}a+\left(  1-\frac{1}{\eta
%	}\right)  b\right)  -E\left(  b\right)  \right]  }  & \geq e^{-\frac{1}{k\eta
%	s}\left[  E\left(  a\right)  -E\left(  b\right)  \right]  }\qquad\forall a\in
%	A\\
%	\frac{1}{\sum_{a\in A}e^{-\frac{1}{ks}\left[  E\left(  \frac{1}{\eta}a+\left(
%	1-\frac{1}{\eta}\right)  b\right)  -E\left(  b\right)  \right]  }}  &
%	\leq\frac{1}{\sum_{a\in A}e^{-\frac{1}{k\eta s}\left[  E\left(  a\right)
%	-E\left(  b\right)  \right]  }}\\
%	p_{s}\left(  b\mid\frac{1}{\eta}A+\left(  1-\frac{1}{\eta}\right)  b\right)
%	& \leq\frac{1}{\sum_{a\in A}e^{-\frac{1}{k\eta s}E\left(  a\right)  +\frac
%	{1}{k\eta s}E\left(  b\right)  }}\\
%	p_{s}\left(  b\mid\frac{1}{\eta}A+\left(  1-\frac{1}{\eta}\right)  b\right)
%	& \leq\frac{e^{-\frac{1}{k\eta s}E\left(  b\right)  }}{\sum_{a\in A}%
%	e^{-\frac{1}{k\eta s}E\left(  a\right)  }}=p_{\eta s}\left(  b\mid A\right)
%	\end{align*}
and $p_{s}\left(  b\mid\left(  1/\eta\right)  A+\left(  1-\left(
1/\eta\right)  \right)  b\right)  \leq p_{\eta s}\left(  b\mid A\right)  $.

\emph{(iii) implies (iv).} Trivial.

\emph{(iv) implies (i).} To prove convexity, it is sufficient to check that,
given any $\alpha\in\left(  0,1\right)  $,%
\begin{equation}
E\left(  \alpha x+\left(  1-\alpha\right)  y\right)  \leq\alpha E\left(
x\right)  +\left(  1-\alpha\right)  E\left(  y\right)  \label{eq:conv}%
\end{equation}
for all $x,y\in X$ such that $E\left(  y\right)  \geq E\left(  x\right)
$.\footnote{In fact, if $E\left(  \bar{x}\right)  >E\left(  \bar{y}\right)  $,
(\ref{eq:conv}) yields, for any $\beta\in\left(  0,1\right)  $, $E\left(
\beta\bar{y}+\left(  1-\beta\right)  \bar{x}\right)  \leq\beta E\left(
\bar{y}\right)  +\left(  1-\beta\right)  E\left(  \bar{x}\right)  $.} Now,
arbitrarily choose $s\in\left(  0,\infty\right)  $. If $E\left(  y\right)
\geq E\left(  x\right)  $, then $y\in\arg\min_{a\in\left\{  x,y\right\}
}p_{s/\alpha}\left(  a\mid\left\{  x,y\right\}  \right)  $, then
(\ref{eq:barba}), with $\eta=1/\alpha$, yields%
\begin{align*}
&  \left.  p_{s}\left(  y\mid\alpha\left\{  x,y\right\}  +\left(
1-\alpha\right)  y\right)  \leq p_{s/\alpha}\left(  y\mid\left\{  x,y\right\}
\right)  \right.  \\
&  \implies r_{s}\left(  y,\alpha x+\left(  1-\alpha\right)  y\right)  \leq
r_{s/\alpha}\left(  y,x\right)  \\
&  \implies\frac{1}{ks}\left[  E\left(  \alpha x+\left(  1-\alpha\right)
y\right)  -E\left(  y\right)  \right]  \leq\frac{\alpha}{ks}\left[  E\left(
x\right)  -E\left(  y\right)  \right]
\end{align*}
for all $\alpha\in\left(  0,1\right)  $, which implies (\ref{eq:conv})%
).\hfill$\blacksquare$\bigskip

\noindent\textbf{Proof of Proposition \ref{prop:bartsuka}} If A.\ref{ax:webo}
is satisfied, by Theorem \ref{thm:Luce_on_steroids} there exist a function
$E:X\rightarrow\mathbb{R}$ and an increasing bijection $\kappa:\left(
0,\infty\right)  \rightarrow\left(  0,\infty\right)  $ such that $p$ is
represented by (\ref{eq:softmax_txt}). Moreover, $p_{\bar{v}}\left(  \bar
{c},\bar{d}\right)  >p_{\bar{v}}\left(  \bar{d},\bar{c}\right)  $ implies
$E\left(  \bar{d}\right)  >E\left(  \bar{c}\right)  $.

For all $a\in A$,%
\begin{align*}
r_{\bar{v}}\left(  a,\bar{c}\right)   &  =\exp\left(  -\frac{1}{\kappa\left(
\bar{v}\right)  }\left[  E\left(  a\right)  -E\left(  \bar{c}\right)  \right]
\right) \\
\frac{1}{\kappa\left(  \bar{v}\right)  }[E\left(  a\right)  -E\left(  \bar
{c}\right)]   &  =-\ln r_{\bar{v}}\left(  a,\bar{c}\right)  =\ln r_{\bar{v}%
}\left(  \bar{c},a\right)
\end{align*}
hence $\ln r_{\bar{v}}\left(  \bar{c},a\right)  =mE\left(  a\right)  +q$, with
$m=1/\kappa\left(  \bar{v}\right)  $ and $q= - E\left(  \bar{c}\right) /\kappa\left(  \bar{v}\right)  $. For
all $t\in\left(  0,\infty\right)  $,%
\[
\frac{\ln r_{\bar{v}}\left(  \bar{c},\bar{d}\right)  }{\ln r_{t}\left(
\bar{c},\bar{d}\right)  }=\frac{-\frac{1}{\kappa\left(  \bar{v}\right)
}\left[  E\left(  \bar{c}\right)  -E\left(  \bar{d}\right)  \right]  }%
{-\frac{1}{\kappa\left(  t\right)  }\left[  E\left(  \bar{c}\right)  -E\left(
\bar{d}\right)  \right]  }=\frac{1}{\kappa\left(  \bar{v}\right)  }%
\kappa\left(  t\right)  =m\kappa\left(  t\right)
\]
Point (i) of Theorem \ref{thm:Luce_on_steroids} implies that
(\ref{eq:softmax_txt}) holds, with $\tilde{E}\left(  \cdot\right)  =mE\left(
\cdot\right)  +q=\ln r_{\bar{v}}\left(  \bar{c},\cdot\right)  \ $%
and\ $\tilde{\kappa}\left(  \cdot\right)  =m\kappa\left(  \cdot\right)  =\ln
r_{\bar{v}}\left(  \bar{c},\bar{d}\right)  /\ln r_{\cdot}\left(  \bar{c}%
,\bar{d}\right)  .$

The converse follows from Theorem \ref{thm:Luce_on_steroids} too: if
representation (\ref{eq:softmax_txt}) holds,\footnote{At the risk of being
pedantic, the sentence \textquotedblleft representation (\ref{eq:softmax_txt})
holds\ for some $\tilde{E}$ and $\tilde{\kappa}$\textquotedblright\ means that
$\tilde{E}:X\rightarrow\mathbb{R}$ is a function, $\tilde{\kappa}:\left(
0,\infty\right)  \rightarrow\left(  0,\infty\right)  $ is an increasing
bijection, and equation (\ref{eq:softmax_txt}) holds for all $\left(
t,a,A\right)  \in\left(  0,\infty\right)  \times X\times\mathcal{A}$.} then
A.\ref{ax:webo} is satisfied.\hspace*{\fill}$\blacksquare$

%\clearpage

\appendix

\section*{Supplementary Material: \\Proof of Proposition
\ref{lem:operationsbis}}

The next Lemma uses the notation of Lemma \ref{lem:operations}.

\begin{lemma}
\label{lem:fantasia}If $p:\left(  0,\infty\right)  \times X\times
\mathcal{A}\rightarrow\mathbb{R}_{+}$ is a random state function that
satisfies A.\ref{ax:pors}, A.\ref{ax:cont}, A.\ref{ax:unif}, and
A.\ref{ax:mono}, then, given any $a,b\in X$:

\begin{itemize}
\item[(i)] $a\succ b$ if and only if $\varphi_{a,b}$ is an increasing
bijection from $\left(  0,\infty\right)  $ to $\left(  1,\infty\right)  $;

\item[(ii)] $a\sim b$ if and only if $\varphi_{a,b}$ is constantly equal to
$1$;

\item[(iii)] $a\prec b$ if and only if $\varphi_{a,b}$ is a decreasing
bijection from $\left(  0,\infty\right)  $ to $\left(  0,1\right)  $.
\end{itemize}

In particular, all the above monotonicity and bijectivity properties are
maintained when $\varphi_{a,b}$ is extended to $\left[  0,\infty\right)  $ by
setting $\varphi_{a,b}\left(  0\right)  =1$.
\end{lemma}

\noindent\textbf{Proof} By the arguments adopted in the proof of Lemma
\ref{lem:operations}, we have that, given any $a,b\in X$, the function
\[%
\begin{array}
[c]{cccc}%
\varphi_{a,b}: & \left(  0,\infty\right)  & \rightarrow & \left(
0,\infty\right) \\
& t & \mapsto & r_{1/t}\left(  a,b\right)
\end{array}
\]
is well defined, and continuous.

\begin{fact}
\label{fct:fiko}If $r_{\tau}\left(  a,b\right)  >1$ for some $\tau\in\left(
0,\infty\right)  $, then
\[%
\begin{array}
[c]{cccc}%
r\left(  a,b\right)  : & \left(  0,\infty\right)  & \rightarrow & \left(
0,\infty\right) \\
& t & \mapsto & r_{t}\left(  a,b\right)
\end{array}
\]
is strictly decreasing and everywhere strictly greater than $1$, that is,
$\varphi_{a,b}$ is strictly increasing and everywhere strictly greater than
$1$.

If $r_{\tau}\left(  a,b\right)  <1$ for some $\tau\in\left(  0,\infty\right)
$, then
\[%
\begin{array}
[c]{cccc}%
r\left(  a,b\right)  : & \left(  0,\infty\right)  & \rightarrow & \left(
0,\infty\right) \\
& t & \mapsto & r_{t}\left(  a,b\right)
\end{array}
\]
is strictly increasing and everywhere strictly smaller than $1$, that is,
$\varphi_{a,b}$ is strictly decreasing and everywhere strictly smaller than
$1$.
\end{fact}

\noindent\textit{Proof} Let $r_{\tau}\left(  a,b\right)  >1$. If $r_{t}\left(
a,b\right)  \leq1$ for some $t>\tau$, by A.\ref{ax:mono} it would follow
$r_{\tau}\left(  a,b\right)  \leq r_{t}\left(  a,b\right)  \leq1$, a
contradiction. Then $r_{t}\left(  a,b\right)  >1$, for all $t\in\left[
\tau,\infty\right)  $. Now, given any $s\in\left(  0,\infty\right)  $, taking
$t\in\left[  \tau,\infty\right)  $ such that $s<t$, by A.\ref{ax:mono} it
follows $r_{s}\left(  a,b\right)  >r_{t}\left(  a,b\right)  >1$. Therefore,
$r_{t}\left(  a,b\right)  >1$, for all $t\in\left(  0,\infty\right)  $. But
then, given any $s<t$ in $\left(  0,\infty\right)  $, since $r_{t}\left(
a,b\right)  >1$, by A.\ref{ax:mono} it follows $r_{s}\left(  a,b\right)
>r_{t}\left(  a,b\right)  $, and
\[%
\begin{array}
[c]{cccc}%
r\left(  a,b\right)  : & \left(  0,\infty\right)  & \rightarrow & \left(
0,\infty\right) \\
& t & \mapsto & r_{t}\left(  a,b\right)
\end{array}
\]
is strictly decreasing, then $\varphi_{a,b}$ is strictly increasing.

Let $r_{\tau}\left(  a,b\right)  <1$, then
\[
r_{\tau}\left(  b,a\right)  =\frac{1}{r_{\tau}\left(  a,b\right)  }>1
\]
hence $r\left(  b,a\right)  $ is strictly decreasing, $r\left(  a,b\right)  $
strictly increasing, $\varphi_{a,b}$ strictly decreasing.\hfill$\square
\medskip$

(i) If $a\succ b$, again by arguments of the proof of Lemma
\ref{lem:operations}, it follows that
\[
\lim_{t\rightarrow
\infty}r_{1/t}\left(  a,b\right) =\lim_{t\rightarrow\infty}\varphi_{a,b}\left(  t\right)  =\infty
\]
then $r_{\tau}\left(  a,b\right)  >1$ for some $\tau\in\left(  0,\infty
\right)  $, and $\varphi_{a,b}$ is strictly increasing. Finally, by A.\ref{ax:mono},
\[
\lim_{t\rightarrow0}\varphi_{a,b}\left(  t\right)  =\lim_{t\rightarrow\infty
}r_{t}\left(  a,b\right)  =1
\]
and so $\varphi_{a,b}$ is an increasing bijection from $\left(  0,\infty
\right)  $ to $\left(  1,\infty\right)  $.

Conversely, if $\varphi_{a,b}$ is an increasing bijection from $\left(
0,\infty\right)  $ to $\left(  1,\infty\right)  $, then
\[
r_{0}\left(  a,b\right)  =\lim_{t\rightarrow0}r_{t}\left(  a,b\right)
=\lim_{t\rightarrow\infty}r_{1/t}\left(  a,b\right)  =\lim_{t\rightarrow
\infty}\varphi_{a,b}\left(  t\right)  =\infty
\]
But then it must be the case that $a\neq b$, and the above limit corresponds
to%
\[
\lim_{t\rightarrow0}\frac{1-p_{t}\left(  b,a\right)  }{p_{t}\left(
b,a\right)  }=\infty
\]
thus $p_{0}\left(  b,a\right)  =0$ and $p_{0}\left(  a,b\right)  =1$. Then, by
definition of $\succsim$, $a\succ b$.

(ii) If $a\sim b$ and $a=b$, then obviously, $\varphi_{a,b}\left(  t\right)
=p_{1/t}\left(  a,b\right)  /p_{1/t}\left(  b,a\right)  =1$, irrespective of
$t\in\left(  0,\infty\right)  $. Else if $a\sim b$ and $a\neq b$, by point (i)
if Lemma \ref{lem:operations} we have that $p_{0}\left(  a,b\right)
=p_{0}\left(  b,a\right)  $, and so%
\begin{equation}
\lim_{t\rightarrow\infty}\varphi_{a,b}\left(  t\right)  =\lim_{t\rightarrow
\infty}\frac{p_{1/t}\left(  a,b\right)  }{p_{1/t}\left(  b,a\right)  }%
=\lim_{t\rightarrow0}\frac{p_{t}\left(  a,b\right)  }{p_{t}\left(  b,a\right)
}=\frac{p_{0}\left(  a,b\right)  }{p_{0}\left(  b,a\right)  }=1 \label{eq:ano}%
\end{equation}
If $\varphi_{a,b}\left(  \bar{t}\right)  >1$\ for some $\bar{t}\in\left(
0,\infty\right)  $, then $r_{\tau}\left(  a,b\right)  >1$ for some $\tau
\in\left(  0,\infty\right)  $ (say, $\tau=1/\bar{t}$), then $\varphi_{a,b}$ is
strictly increasing, which contradicts (\ref{eq:ano}), because it implies
$\lim_{t\rightarrow\infty}\varphi_{a,b}\left(  t\right)  \geq\varphi
_{a,b}\left(  \bar{t}\right)  >1$. If $\varphi_{a,b}\left(  \bar{t}\right)
<1$\ for some $\bar{t}\in\left(  0,\infty\right)  $, then $r_{\tau}\left(
a,b\right)  <1$ for some $\tau\in\left(  0,\infty\right)  $ (say, $\tau
=1/\bar{t}$), then $\varphi_{a,b}$ is strictly decreasing, which contradicts
(\ref{eq:ano}). Therefore, $\varphi_{a,b}\left(  t\right)  =1$, irrespective
of $t\in\left(  0,\infty\right)  $.

Conversely, if $\varphi_{a,b}\left(  t\right)  \equiv1$, then $r_{t}\left(
a,b\right)  \equiv1$, hence $p_{t}\left(  a,b\right)  \equiv p_{t}\left(
b,a\right)  $, and $p_{0}\left(  a,b\right)  =p_{0}\left(  b,a\right)  $, thus
$a\sim b$.

(iii) $a\prec b$ iff $b\succ a$ iff $\varphi_{b,a}$ is an increasing bijection
from $\left(  0,\infty\right)  $ to $\left(  1,\infty\right)  $ iff
$\varphi_{a,b}=1/\varphi_{b,a}$ is a decreasing bijection from $\left(
0,\infty\right)  $ to $\left(  0,1\right)  $. \hfill$\blacksquare\bigskip$

By the previous arguments, and since, by Lemma \ref{lem:operations}, $\succsim$ is a trichotomy, we have the following:

\begin{corollary}
\label{lem:fanta}If $p:\left(  0,\infty\right)  \times X\times\mathcal{A}%
\rightarrow\mathbb{R}_{+}$ is a random state function that satisfies
A.\ref{ax:pors}, A.\ref{ax:cont}, A.\ref{ax:unif}, and A.\ref{ax:mono}, then,
given any $a,b\in X$:

\begin{itemize}
\item[(i)] $a\succ b$ if and only if $r_{t}\left(  a,b\right)  >1$ for some/all $t\in\left(  0,\infty\right)  $;

\item[(ii)] $a\sim b$ if and only if $r_{t}\left(  a,b\right)  =1$ for some/all $t\in\left(  0,\infty\right)  $;

\item[(iii)] $a\prec b$ if and only if $r_{t}\left(  a,b\right)  <1$ for some/all $t\in\left(  0,\infty\right)  $.
\end{itemize}
\end{corollary}

\noindent\textbf{Proof of Proposition \ref{lem:operationsbis}} Assume $p$ is
not uniform (the uniform case is left to the reader). If $p$ satisfies
A.\ref{ax:cons} and A.\ref{ax:webo}, then using the representation provided by
Theorem \ref{thm:Luce_on_steroids}, it is routine to show that it satisfies
A.\ref{ax:mono} and A.\ref{ax:conc}. We only prove the converse.

As to A.\ref{ax:cons}, let $\left(  t,a,b\right)  \in\left(  0,\infty\right)
\times X^{2}$ be such that $p_{t}\left(  a,b\right)  >p_{t}\left(  b,a\right)
$. Then $a\neq b$ and $r_{t}\left(  a,b\right)  >1$, by the previous results,
$\varphi_{a,b}\left(  t\right)  =r_{1/t}\left(  a,b\right)  $ is an increasing
bijection from $\left(  0,\infty\right)  $ to $\left(  1,\infty\right)  $,
then%
\[
\lim_{s\rightarrow0}\frac{p_{s}\left(  a,b\right)  }{1-p_{s}\left(
a,b\right)  }=\lim_{s\rightarrow0}r_{s}\left(  a,b\right)  =\lim
_{t\rightarrow\infty}\varphi_{a,b}\left(  t\right)  =\infty
\]
thus $p_{0}\left(  a,b\right)  =1>0=p_{0}\left(  b,a\right)  $. As wanted.

As to A.\ref{ax:webo}. Given any $a,b\in X$, set $\varphi_{a,b}\left(
0\right)  =1$ as in Lemma \ref{lem:fantasia}. Denote $w_{t}\left(  a,b\right)
=\ln\varphi_{a,b}\left(  t\right)  $, for all $\left(  t,a,b\right)
\in\left[  0,\infty\right)  \times X^{2}$. Arbitrarily choose $\hat{a}%
\succ\hat{b}\in X$, so that $\varphi_{\hat{a},\hat{b}}:\left[  0,\infty
\right)  \rightarrow\left[  1,\infty\right)  $ is an increasing bijection, and
notice that the function%
\begin{equation}
f\left(  t\right)  =\ln\varphi_{\hat{a},\hat{b}}\left(  t\right)  =w_{t}%
(\hat{a},\hat{b})\qquad\forall t\in\left[  0,\infty\right)  \label{eq:bij}%
\end{equation}
is an increasing bijection onto $\left[  0,\infty\right)  $, so $f\left(
0\right)  =0$ and $f_{|\left(  0,\infty\right)  }$ is an increasing bijection
onto $\left(  0,\infty\right)  $. The next steps verify that $p$ satisfies
A.\ref{ax:webo} with respect to $\oplus_{f}$.

Note that, given any $t,s\in\left(  0,\infty\right)  $, we have%
\begin{align}
w\underset{\tau}{\underbrace{_{f^{-1}\left(  f\left(  t\right)  +f\left(
s\right)  \right)  }}}(\hat{a},\hat{b})  &  =f\left(  \underset{\tau
}{\underbrace{f^{-1}\left(  f\left(  t\right)  +f\left(  s\right)  \right)  }%
}\right) \label{eq:ICAcause}\\
&  =f\left(  t\right)  +f\left(  s\right)  =w_{t}(\hat{a},\hat{b})+w_{s}%
(\hat{a},\hat{b})
\end{align}
Next we show that (\ref{eq:ICAcause}) and A.\ref{ax:conc} imply%
\begin{equation}
w_{f^{-1}\left(  f\left(  t\right)  +f\left(  s\right)  \right)  }\left(
a,b\right)  =w_{t}\left(  a,b\right)  +w_{s}\left(  a,b\right)
\label{eq:ICAeffect}%
\end{equation}
for all $a,b\in X$ and all $t,s\in\left(  0,\infty\right)  $. Given{\ any
$c,d,x,y\in X$ and }any $s,t,\tau\in\left(  0,\infty\right)  $ such that
$w_{\tau}\left(  c,d\right)  >0$ and $w_{\tau}\left(  x,y\right)  >0$, we have
$r_{1/\tau}\left(  c,d\right)  =e^{w_{\tau}\left(  c,d\right)  }>1$ and
$r_{1/\tau}\left(  x,y\right)  =e^{w_{\tau}\left(  x,y\right)  }>1$, hence, by
A.\ref{ax:conc},
\begin{align*}
r_{1/\tau}\left(  c,d\right)   &  >r_{1/t}\left(  c,d\right)  r_{1/s}\left(
c,d\right) \\
&  \iff r_{1/\tau}\left(  x,y\right)  >r_{1/t}\left(  x,y\right)
r_{1/s}\left(  x,y\right) \\
w_{t}\left(  c,d\right)   &  >w_{t}\left(  c,d\right)  +w_{s}\left(
c,d\right) \\
&  \iff w_{\tau}\left(  x,y\right)  >w_{t}\left(  x,y\right)  +w_{s}\left(
x,y\right)
\end{align*}
(the roles of $\left(  c,d\right)  $ and $\left(  x,y\right)  $ are symmetric
in the axiom). By Corollary \ref{lem:fanta}, if $w_{\hat{t}}\left(
c,d\right)  >0$ and $w_{\hat{s}}\left(  x,y\right)  >0$ for some $\hat{t}%
,\hat{s}\in\left(  0,\infty\right)  $, then $w_{\tau}\left(  c,d\right)  >0$
and $w_{\tau}\left(  x,y\right)  >0$ for all $\tau\in\left(  0,\infty\right)
$. Therefore, given any{\ $c,d,x,y\in X$, }if $w_{\hat{t}}\left(  c,d\right)
>0$ and $w_{\hat{s}}\left(  x,y\right)  >0$ for some $\hat{t},\hat{s}%
\in\left(  0,\infty\right)  $, then, given any {$s,t,\tau\in\left(
0,\infty\right)  $}, it follows{%
\begin{align}
w_{\tau}\left(  c,d\right)   &  \leq w_{t}\left(  c,d\right)  +w_{s}\left(
c,d\right) \label{eq:acci}\\
&  \iff w_{\tau}\left(  x,y\right)  \leq w_{t}\left(  x,y\right)
+w_{s}\left(  x,y\right)
\end{align}
}Moreover, as we argued for (\ref{eq:bij}), since $c\succ d$ and $x\succ y$,
the functions $h\left(  t\right)  =w_{t}(c,d)$ and $g\left(  t\right)
=w_{t}(x,y)$ are increasing bijections from $\left(  0,\infty\right)  $ to
$\left(  0,\infty\right)  $ and (\ref{eq:acci}) implies
\[
\tau\leq h^{-1}\left(  h\left(  t\right)  +h\left(  s\right)  \right)
\iff\tau\leq g^{-1}\left(  g\left(  t\right)  +g\left(  s\right)  \right)
\]
for all $s,t,\tau\in\left(  0,\infty\right)  $. But then, $h^{-1}\left(
h\left(  t\right)  +h\left(  s\right)  \right)  =g^{-1}\left(  g\left(
t\right)  +g\left(  s\right)  \right)  $ for all $s,t\in\left(  0,\infty
\right)  $. Hence, for all $s,t,\tau\in\left(  0,\infty\right)  $,%
\[
\tau=h^{-1}\left(  h\left(  t\right)  +h\left(  s\right)  \right)  \iff
\tau=g^{-1}\left(  g\left(  t\right)  +g\left(  s\right)  \right)
\]
that is, $h\left(  \tau\right)  =h\left(  t\right)  +h\left(  s\right)  \iff
g\left(  \tau\right)  =g\left(  t\right)  +g\left(  s\right)  $.

Therefore:

\begin{itemize}
\item if $w_{\hat{t}}\left(  c,d\right)  >0$ and $w_{\hat{s}}\left(
x,y\right)  >0$ for some $\hat{t},\hat{s}\in\left(  0,\infty\right)  $, then,
given any {$s,t,\tau\in\left(  0,\infty\right)  $}, it holds
\begin{align*}
w_{\tau}\left(  c,d\right)   &  =w_{t}\left(  c,d\right)  +w_{s}\left(
c,d\right) \\
&  \iff w_{\tau}\left(  x,y\right)  =w_{t}\left(  x,y\right)  +w_{s}\left(
x,y\right)
\end{align*}

\item if $w_{\hat{t}}\left(  c,d\right)  >0$ and $w_{\hat{s}}\left(
x,y\right)  <0$ for some $\hat{t},\hat{s}\in\left(  0,\infty\right)  $, then,
$w_{\hat{s}}\left(  y,x\right)  >0$ and, given any {$s,t,\tau\in\left(
0,\infty\right)  $}, it holds
\begin{align*}
w_{\tau}\left(  c,d\right)   &  =w_{t}\left(  c,d\right)  +w_{s}\left(
c,d\right) \\
&  \iff w_{\tau}\left(  y,x\right)  =w_{t}\left(  y,x\right)  +w_{s}\left(
y,x\right) \\
&  \iff-w_{\tau}\left(  y,x\right)  =-w_{t}\left(  y,x\right)  -w_{s}\left(
y,x\right) \\
&  \iff w_{\tau}\left(  x,y\right)  =w_{t}\left(  x,y\right)  +w_{s}\left(
x,y\right)
\end{align*}

\item if $w_{\hat{t}}\left(  c,d\right)  >0$ and $w_{\hat{s}}\left(
x,y\right)  =0$ for some $\hat{t},\hat{s}\in\left(  0,\infty\right)  $, then,
$\varphi_{x,y}$ is constantly equal to $1$, and $w_{\tau}\left(  x,y\right)
=w_{t}\left(  x,y\right)  =w_{s}\left(  x,y\right)  =0$, for all {$s,t,\tau
\in\left(  0,\infty\right)  $}, thus, given any {$s,t,\tau\in\left(
0,\infty\right)  $}, it holds
\begin{align*}
w_{\tau}\left(  c,d\right)   &  =w_{t}\left(  c,d\right)  +w_{s}\left(
c,d\right) \\
&  \implies w_{\tau}\left(  x,y\right)  =w_{t}\left(  x,y\right)
+w_{s}\left(  x,y\right)
\end{align*}

\end{itemize}

Summing up, since $\hat{a}\succ\hat{b}$, then, given any {$s,t,\tau\in\left(
0,\infty\right)  $,}
\begin{align}
w_{\tau}(\hat{a},\hat{b})  &  =w_{t}(\hat{a},\hat{b})+w_{s}(\hat{a},\hat
{b})\label{Original_ICA}\\
&  \implies w_{\tau}\left(  x,y\right)  =w_{t}\left(  x,y\right)
+w_{s}\left(  x,y\right)
\end{align}
for all $x,y\in X$. Now by (\ref{eq:ICAcause})
\[
w\underset{\tau}{\underbrace{_{f^{-1}\left(  f\left(  t\right)  +f\left(
s\right)  \right)  }}}(\hat{a},\hat{b})=w_{t}(\hat{a},\hat{b})+w_{s}(\hat
{a},\hat{b})\qquad\forall t,s\in\left(  0,\infty\right)
\]
and so (\ref{Original_ICA}) implies
\[
w\underset{\tau}{\underbrace{_{f^{-1}\left(  f\left(  t\right)  +f\left(
s\right)  \right)  }}}\left(  x,y\right)  =w_{t}\left(  x,y\right)
+w_{s}\left(  x,y\right)  \qquad\forall t,s\in\left(  0,\infty\right)
\]
and for all $x,y\in X$. Finally, for all $x,y\in X$ and all $t,s\in\left(
0,\infty\right)  $%
\begin{align*}
r_{\frac{1}{t\oplus_{f}s}}\left(  x,y\right)   &  =\varphi_{x,y}\left(
t\oplus_{f}s\right)  =\varphi_{x,y}\left(  f^{-1}\left(  f\left(  t\right)
+f\left(  s\right)  \right)  \right) \\
&  =e^{w_{f^{-1}\left(  f\left(  t\right)  +f\left(  s\right)  \right)
}\left(  x,y\right)  }=e^{w_{t}\left(  x,y\right)  }e^{w_{s}\left(
x,y\right)  }\\
&  =\varphi_{x,y}\left(  t\right)  \varphi_{x,y}\left(  s\right)  =r_{\frac
{1}{t}}\left(  x,y\right)  r_{\frac{1}{s}}\left(  x,y\right)
\end{align*}
and A.\ref{ax:webo} holds.\hfill$\blacksquare$

\end{document}